\documentclass[12pt]{article}
\usepackage{amssymb}
\usepackage{amsmath}
\usepackage{amsthm}

\newtheorem{theorem}{Theorem}[section]

\newtheorem{proposition}[theorem]{Proposition}

\theoremstyle{definition}

\theoremstyle{remark}

\numberwithin{equation}{section}

\usepackage[mathcal]{eucal}
\usepackage{hyperref}

\newcommand{\ve}{\varepsilon}
\newcommand{\la}{\langle}
\newcommand{\ra}{\rangle}
\newcommand{\tos}{\rightrightarrows}

\newcommand{\cpl}{\Phi}
\newcommand{\J}{\mathcal{J}}
\newcommand{\ch}{\mathcal{H}}

\title{Fixed points in the
  family of convex representations of a maximal monotone operator\\
{\small published on:
  \href{http://www.ams.org/proc/2003-131-12/S0002-9939-03-07083-7/home.html}{
         {\it Proc.\ Amer.\ Math.\ Soc.\ } {\bf 131} (2003) 3851--3859.
         }
}
}

\author{{\it B. F. Svaiter}\thanks{Partially supported by 
                             CNPq Grant 301200/93-9(RN)
                             and by PRONEX--Optimization.}
         \\ IMPA Instituto de Matem\' atica Pura e Aplicada
         \\ Estrada Dona Castorina 110
         \\ Rio de Janeiro--RJ
         \\ CEP 22460-320 Brazil
         \\ email: benar@impa.br}

\date{7 August 2002}

\begin{document}

\maketitle

\begin{abstract}
Any  maximal monotone operator can be characterized by a convex
function.  The family of such convex functions is invariant under a
transformation connected with the Fenchel-Legendre conjugation.  We
prove that there exist a convex representation of the operator which
is a fixed point of this conjugation.
\\
\\
2000 Mathematics Subject Classification: 47H05
\\
\\
{ \sc keywords:} maximal monotone operators, conjugation, convex functions
\end{abstract}

\section{Introduction}
Let $X$ be a real Banach space and $X^*$ its dual.  It is usual to
identify a point to set operator $T:X\tos X^*$ with its graph, $\{
(x,x^*)\in X\times  X^*\,|\, x^*\in T(x)\}$.  We will use the notation $\la
x,x^*\ra$ for the duality product $x^*(x)$ of $x\in X$, $x^*\in X^*$.

An operator $T:X\tos X^* $ is \emph{monotone} if
\[ (x,x^*),(y, y^*)\in T\Rightarrow \la x-y, x^*-y^*\ra\geq 0,\]
and is 
is \emph{maximal  monotone} if it is monotone and
\[ \forall (y,y^*)\in T,\la x-y, x^*-y^*\ra\geq 0\Rightarrow
(x,x^*)\in T.\]

Krauss~\cite{krauss} managed to represent maximal monotone operators by
subdifferentials of saddle functions on $X\times X$.  After that,
Fitzpatrick~\cite{fitzpatrick:1988} proved that maximal monotone
operators can be represented by convex functions on $X\times X^*$.
Latter on, Simons~\cite{simons.book} studied maximal monotone
operators using a min-max approach.  Recently, the convex representation
of maximal monotone operators was rediscovered by Burachik and
Svaiter~\cite{burachik-svaiter:2002} and Martinez-Legaz and
Th\'era~\cite{martinez-thera:2001}.  In \cite{burachik-svaiter:2002},
some results on enlargements are used to perform a systematic study of
the family of convex functions which represents a given maximal
monotone operator.  Here we are concerned with this kind of
representation.

Given $f:X\to\overline{\mathbb{R}}$, the \emph{Fenchel-Legendre}
conjugate of $f$ is $f^*:X^*\to\overline{\mathbb{R}}$,
\[ f^*(x^*):=\sup_{x\in x} \la x,x^*\ra- f(x).\]
The \emph{subdifferential} of $f$ is the operator $\partial f:X\tos
X^*$,
\[ \partial f(x):=\{ x^*\in X^*\,|\, f(y)\geq f(x)+\la y-x,x^*\ra, \,
\forall y\in X\}.\]
If $f$ is convex, lower semicontinuous and proper, then $\partial f$
is maximal monotone~\cite{rockafellar:1970mmsub}.  From the previous
definitions, we have the \emph{Fenchel--Young inequality:} for all
$x\in X$, $x^*\in X^*$
\[f(x)+f^*(x^*)\geq \la x, x^*\ra\, ,\:\:
f(x)+f^*(x^*)= \la x, x^*\ra \iff x^*\in\partial f(x).
\]
So, defining $h_{\mathrm{FY}}:X\times X^*\to\overline{\mathbb{R}}$,
\begin{equation}
  \label{eq:new1}
   h_{\mathrm{FY}}(x,x^*):=f(x)+f^*(x^*),
\end{equation}
we observe that this function fully characterizes $\partial f$.
Assume that $f$ is convex, lower semicontinuous and proper. In this
case, $\partial f$ is maximal monotone.  Moreover, if we use the
canonical injection of $X$ in to $X^{**}$, then $f^{**}(x)=f(x)$ for all
$x\in X$. Hence, for all $(x,x^*)\in X\times X^*$
\[  (h_{\mathrm{FY}})^*(x,x^*)=h_{\mathrm{FY}}(x,x^*).\]
Our aim it to prove that any maximal monotone operator has a convex
representation with a similar property.

From now on, $T:X\tos X^*$ is a maximal monotone operator.  Define, as
in \cite{fitzpatrick:1988}, $\mathcal{H}(T)$ to be the family of
convex lower semi continuous functions $h:X\times X^* \to
\overline{\mathbb{R}}$ such that
 \begin{equation}
  \label{eq:defh}
  \begin{array}{rl}
  \forall (x,x^*)\in X\times X^*,&\;  h(x,x^*)\geq \la x, x^*\ra,\\
 &  (x,x^*)\in T\Rightarrow h(x,x^*)=\la x,x^*\ra.
\end{array}
\end{equation}
This family is nonempty \cite{fitzpatrick:1988}.
Moreover, for any $h\in\mathcal{H}(T)$, \( h(x,x^*)=\la x,x^*\ra\) if
and only if $(x,x^*)\in T$~\cite{burachik-svaiter:2002}.  Hence, any
element of $\mathcal{H}(T)$ fully characterizes, or represents $T$.
Since the $\sup$ of convex lower semicontinuous function is also
convex and lower semicontinuous, using also 
\eqref{eq:defh} we
conclude that $\sup$ of any (nonempty) subfamily of $\ch(T)$ is still
in $\ch(T)$.

The dual of $X\times X^{*}$ is $X^*\times X^{**}$. So, for
$(x,x^*)\in X\times X^{*}$, $(y^*,y^{**})\in X^*\times X^{**}$,
\[ \la (x,x^*)\,,\, (y^*,y^{**})\ra=\la x,y^*\ra+\la x^*, y^{**}\ra.\]
Given an function $h:X\times X^*\to\overline{\mathbb{R}}$, define
$\J h:X\times X^*\to\overline{\mathbb{R}}$,
\begin{equation}
  \label{eq:defj}
\J h(x,x^*):=h^*(x^*,x),
\end{equation}
where $h^*$ stands for the Fenchel-Legendre conjugate of $h$ and the
canonical inclusion of $X$ in $X^{**}$ is being used.
Equivalently,
\begin{equation}
  \label{eq:defj2}
   \J h(x,x^*)=\sup_{(y,y^*)\in X\times X^*}
   \la x,y^*\ra+\la y,x^*\ra-h(y,y^*).
\end{equation}
Trivially, $\J$ inverts the natural order of functions, i.e., if
$h\geq h'$ then $\J h'\geq \J h$.
The family
$\mathcal{H}(T)$ is invariant under the application $\J
$~\cite{burachik-svaiter:2002}. The aim of this paper is to prove that
there exist an element $h\in\mathcal{H}(T)$ such that $\J h=h$.

The application  $\J $ can be studied in the framework of \emph{generalized
  conjugation}~\cite[Ch.\ 11, Sec.\ L]{rock-wets}.
With this aim,
define
\[\begin{array}{l}
\cpl:(X\times X^*)\times (X\times X^*):\to\mathbb{R},\\
\cpl((x,x^*),(y,y^*)):=\la x,y^*\ra+\la y,x^*\ra.
\end{array}
\]
Given $h:X\times X^*\to \overline{\mathbb{R}}$, let $h^\cpl$ be the
conjugate of $h$ with respect to the coupling function $\cpl$,
\begin{equation}
  \label{eq:conj1}
  h^\cpl(x,x^*):=\sup_{(y,y^*)\in X\times X^*} \cpl((x,x^*),(y,y^*))-h(y,y^*).
\end{equation}
Now we have
\[ \J h =h^\cpl,\]
and, in particular
\begin{equation}
  \label{eq:conj2}
  h\geq h^{\cpl\cpl}= \J ^2 h.
\end{equation}

\section{Proof of the Main Theorem}

Define as in \cite{burachik-svaiter:2002},
$\sigma_T:X\times X^*\to\overline{\mathbb{R}}$,
\[\sigma_T:=\sup_{h\in \ch(T)}\; h.\]
Since $\ch(T)$ is ``closed'' under the $\sup$ operation, we conclude
that $\sigma_T$ is the biggest element of $\ch(T)$. Combining this fact
with the inclusion $\J \sigma_T\in \ch(T)$ we conclude that
\[ \sigma_T\geq \J \sigma_T.\]
For a more detailed discussion on  $\sigma_T$, we refer the reader to
\cite[eq.\ (35)]{burachik-svaiter:2002}.  The above inequality will be, in some
sense our departure point.
Define now
\[ \ch_a(T):=\{ h\in\ch(T)\,|\, h\geq \J h \}.\]
The family $\ch_a(T)$ is connected with a family of enlargements of
$T$ which shares with the $\varepsilon$-subdifferential a special
property (see~\cite{{burachik-svaiter:2002}}).  We already know that
$\sigma_T\in \ch_a(T)$. Latter on, we will use the following
construction of elements in this set.
\begin{proposition}
\label{pr:gen}
Take $h\in\ch(T)$ and define
\[ \hat h=\max\;  h,\, \J h.\]
Then $\hat h\in \ch_a(T)$.
\end{proposition}
\begin{proof}
Since $h$ and $\J h$ are in $\ch(T)$, $\hat h\in
\ch(T)$.
By definition,
\[ \hat h\geq h,\; \hat h\geq \J h.\]
Applying $\J$ on these inequalities and using (\ref{eq:conj2}) for
majorizing $\J^2 h$ we obtain
\[ \J h\geq \J \hat h,\; h\geq \J \hat h.\]
Hence, $\hat h\geq \J\hat h$.
\end{proof}

For $h\in \ch(T)$ define
\[ L(h):=\{ g\in\mathcal{H}(T)\,|\, h\geq g\geq \J g\}.\]
The operator $\J $ inverts the order. Therefore,
$L(h)\neq\emptyset$ if and only if $h\geq\J h$, i.e., $h\in\ch_a(T)$.
We already know that $L(\sigma_T)\neq\emptyset$.
\begin{proposition}
For any  $h\in\mathcal{H}_a(T)$, the family $L(h)$ has a
minimal element.
\end{proposition}
\begin{proof}
We shall use Zorn Lemma.  Let $\mathcal C\subseteq L(h)$
 be
a (nonempty) chain, that is, $\mathcal{C}$ is totally ordered.
Take  $h'\in\mathcal{C}$. For any $h''\in\mathcal{C}$,
$h'\geq h''$ or $h''\geq h'$. In the first case we have
$ h'\geq h''\geq \J h''$,
and in the second case,
$ h'\ge \J h'\geq  \J h''$. Therefore,
\begin{equation}
  \label{eq:pr1.r1}
  h'\geq \J h'',\;\forall h',h''\in \mathcal{C}.
\end{equation}
Define now
\begin{equation}
  \label{eq:pr1.dg}
  \hat g =\sup_{h'\in \mathcal C} \J h'.
\end{equation}
Since $\ch(T)$ is invariant under $\J $ and also closed
with respect to the $\sup$, we have $\hat g\in \mathcal{H}(T)$.
From (\ref{eq:pr1.r1}), (\ref{eq:pr1.dg}) it follows that
\[ h'\geq \hat g\geq \J h',\;\forall h'\in\mathcal{C}.\]
Applying $\J $ on the above inequalities, and using also
(\ref{eq:conj2}), we conclude that,
\begin{equation}
  \label{eq:pr1.r3}
  h'\geq \J \hat g\geq \J h',\;
  \forall h'\in \mathcal{C}.
\end{equation}
Since $\hat g\in\mathcal{H}(T)$, $\J \hat g\in\mathcal{H}(T)$.  Taking
the $\sup$ on $h'\in\mathcal{C}$, in the right had side of the last
inequality, we get
\[ \J \hat g\geq\hat g.\]
Applying $\J $, again, we obtain
\[ \J \hat g \geq \J ( \J \hat g).\]
Take some $h' \in\mathcal{C}$. By the definition of $L(h)$ and
(\ref{eq:pr1.r3}), we conclude that $h\geq h'\geq \J \hat g$. Hence
$\J \hat g$ belongs to $L(h)$ and is a lower bound for any element of
$\mathcal{C}$.
Now we apply Zorn Lemma to conclude that $L(h)$ has a minimal
element.
\end{proof}

The minimal elements of $L(h)$ (for $h\in\ch_a(T)$) are the natural
candidates for being fixed points of $\J$.  First we will show that
they are fixed points of $\J^2$.  Observe that, since $\J$
inverts the order of functions, $\J^2$ preserves it, i.e., if $h\geq
h'$ then $\J^2 h\geq \J^2 h'$.  Moreover, $\J^2$ maps $\ch(T)$ in itself.

\begin{proposition}\label{pr:fj2}
  Take $h\in\ch_a(T)$ and let $h_0$ be a minimal element of $L(h)$.
  Then $\J^2 h_0=h_0$.
\end{proposition}
\begin{proof}
First observe that $\J^2 h_0\in\ch(T)$.
By assumption, $h_0\geq \J h_0$.  Applying $\J^2$ in
this inequality we get
\[\J^2 h_0 \geq \J^2( \J h_0)=\J(\J^2 h_0).\]
Since $h\geq h_0$ and, by (\ref{eq:conj2}) $h_0\geq \J^2 h_0$, we
conclude that $h\geq \J^2 h_0\geq \J(\J^2 h_0)$.  Hence $\J^2 h_0\in
L(h)$.
Using again the inequality $h_0\geq \J^2 h_0$ and the minimality of
$h_0$, the conclusion follows.
\end{proof}

\begin{theorem}\label{th:main}
  Take $h\in\mathcal{H}(T)$ such that $h\geq\J h$.  Then
  $h_0\in L(h)$ is minimal (on $L(h)$) if and only if
  $h_0=\J h_0$.
\end{theorem}
\begin{proof}
Assume first that $h_0=\J h_0$. If $h'\in L(h)$ and
\[ h_0\geq h',\]
then, applying $\J $ on this inequality and using the definition
of $L(h)$ we conclude that
\[ h'\geq \J h'\geq \J h_0=h_0.\]
Combining the above inequalities we obtain $h'=h_0$. Hence $h_0$ is
minimal on $L(h)$.

Assume now that $h_0$ is minimal on $L(h)$.
By the definition of $L(h)$, $h_0\geq \J h_0$. Suppose that for some
$(x_0,x_0 ^*)$,
\begin{equation}
  \label{eq:abs}
  h_0(x_0,x_0 ^*)>\J h_0(x_0,x_0 ^*).
\end{equation}
We shall prove that this assumption is contradictory.  By Proposition
\ref{pr:fj2}, $h_0=\J (\J h_0)$. Hence, the above inequality can be
expressed as
\[ \J(\J h_0)(x_0,x_0 ^*)> \J h_0 (x_0, x_0 ^*),\]
or equivalently
\[ \sup_{(y,y^*)\in X\times X^*}
\la y,x_0 ^*\ra+ \la x_0, y^*\ra-\J h_0(y,y^*)> \J h_0 (x_0, x_0 ^*).\]
Therefore, there exists some $(y_0, y_0 ^*)\in X\times X^*$ such that
\begin{equation}
  \label{eq:in1}
 \la y_0,x_0 ^*\ra+\la x_0 ,y_0 ^*\ra-\J h_0 (y_0 ,y_0 ^*)>
 \J h_0 (x_0 ,x_0 ^*).
\end{equation}
In particular, $\J h_0 (y_0,y_0 ^*),\J h_0  (x_0 ,x_0 ^*) \in
\mathbb{R}$. Interchanging
$\J h_0(y_0 ,y_0 ^*)$ with $\J h_0( x_0 ,x_0 ^*)$ we get
\[ \la y_0 ,x_0 ^*\ra+\la x_0 ,y_0 ^*\ra-\J h_0 (x_0,x_0 ^*)>
\J h_0 (y_0 ,y_0 ^*).\]
Therefore, using also (\ref{eq:defj2}), we get
$\J(\J h_0 (y_0,y_0 ^*))> \J h_0(y_0, y_0 ^*)$.
Using again the equality $\J^2 h_0 =h_0$ we conclude that
\begin{equation}
  \label{eq:in2}
   h_0(y_0,y_0 ^*) 
 >\J h_0 (y_0,y_0 ^*).
\end{equation}

Define  $\gamma:X\times X^*\to\mathbb{R}$,
$g:X\times X^*\to\overline{\mathbb{R}}$,
\begin{eqnarray}
  \label{eq:in3a}
    \gamma(x,x^*)&:=&\la x,y_0 ^*\ra+
       \la y_0 , x^*\ra-\J h_0(y_0,y_0 ^*),\\[.2em]
   \label{eq:in3b}
   g&:=&\max\; \gamma,\, \J h_0.
\end{eqnarray}
By (\ref{eq:defj2}), $h_0\geq\gamma$.  Since $h_0\in L(h)$, $h_0\geq
\J h_0$. Therefore,
\[ h_0\geq g \geq \J h_0.\]
We claim that $g\in\ch(T)$. Indeed, $g$ is a lower semicontinuous
convex function.  Moreover,
since
$h_0,\J h_0\in \ch (T)$, it follows from (\ref{eq:defh}) and the above
inequalities that $g\in\ch(T)$.  Now apply $\J$ to the above inequality
to conclude that
\[ h_0\geq \J g\geq \J h_0.\]
Therefore, defining
\begin{equation}
  \label{eq:defgh}
  \hat g=\max\; g,\, \J g,
\end{equation}
we have $h>h_0\geq \hat g$. By Proposition \ref{pr:gen},
$\hat g\in \ch(T)$ and $\hat g\geq
\J\hat g$.  Combining these results with the minimality of $h_0$, it follows
that $\hat g=h_0$.  In particular,
\begin{equation}
  \label{eq:eq5}
   \hat g(y_0,y_0 ^*) = h_0(y_0 ,y_0 ^*).
\end{equation}
To end the prove we shall evaluate $\hat g(y_0,y_0 ^*)$. Using
(\ref{eq:in3a}) 
we obtain
\[ \gamma(y_0 ,y_0 ^*)=2\la y_0 ,y_0 ^*\ra-\J h_0(y_0 ,y_0 ^*).
\]
Since $\J h_0\in \ch(T)$, $\J h_0(y_0, y_0 ^*)\geq \la y_0, y_0
^*\ra$.
Hence, $\gamma(y_0 ,y_0 ^*)\leq \la y_0, y_0 ^*\ra$ and by (\ref{eq:in3b})
\begin{equation}
  \label{eq:in7}
   g(y_0,y_0 ^*)=\J h_0(y,y^*).
\end{equation}
Using again the inequality $g\geq \gamma$, we have
\[ \J \gamma (y_0,y_0 ^*)\geq \J g(y_0 ,y_0 ^*).\]
Direct calculation yields $\J \gamma (y_0,y_0 ^*)=\J h_0 (y,y^*)$. Therefore
\begin{equation}
  \label{eq:in9}
  \J h_0 (y_0,y_0 ^*)\geq \J g(y_0,y_0 ^*).
\end{equation}
Combining (\ref{eq:in7}), (\ref{eq:in9}) and (\ref{eq:defgh}) we obtain
\[ \hat g (y_0,y_0 ^*)=\J h_0 (y_0 ,y_0 ^*).\]
This equality, together with (\ref{eq:eq5}) yields $h_0(y_0,y_0 ^*)=
\J h_0 (y_0, y_0 ^*)$, in contradiction with (\ref{eq:in2}).
Therefore, $h_0(x,x^*)=\J h_0(x,x^*)$ for all $(x,x^*)$.
\end{proof}

Since $\sigma_T\in \ch_a(T)$, $L(\sigma_T)\neq\emptyset$ and there
exist some $h\in L(\sigma_T)$ such that $\J h=h$. (Indeed
$L(\sigma_T)=\ch_a(T)$.)

\section{Application}

Let $f:X\tos X^*$ be a proper lower semicontinuous convex function. We
already know that $\partial f$ is maximal monotone.  Define, for
$\ve\geq 0$,
\[ \partial_\ve f(x):=\{ x^*\in X^*\,|\, f(y)\geq f(x)+\la y-x,x^*\ra-\ve, \,
\forall y\in X\}.\]
Note that $\partial_{0} f=\partial f$. We also have
\begin{eqnarray}
  \label{eq:a1}
  &&\partial f(x)\subseteq \partial_\ve f(x),\, \forall x\in X,\ve\geq 0,\\
\label{eq:a2}
  &&0\leq \ve_1\leq\ve_2\Rightarrow \partial_{\ve_1} f(x)
  \subseteq \partial_{\ve_2} f(x),\, \forall x\in X
\end{eqnarray}
Property (\ref{eq:a1}) tells that $\partial_\ve f$ \emph{enlarges}
$\partial f$. Property (\ref{eq:a2}) shows that $\partial_\ve f$ is
nondecreasing  (or increasing) in $\ve$.   The operator $\partial_\ve
f$ has been introduced in \cite{brock:65}, and since that, it has had may
theoretical and algorithmic applications
\cite{bertsekas-mitter-1973,nurminski-1986,hiriart-lemarechal-1993,
kiwiel-1990,schramm-zowe-1992,lemarechal-nemirovskii-nesterov-1995,
bgls-1997}.

 Since $\partial f$
is maximal monotone, the enlarged operator $\partial_\ve f$ loses
monotonicity in general. Even though, we have
\begin{equation}
  \label{eq:a3}
  x^*\in\partial_\ve f(x)\Rightarrow \la x-y,x^*-y^*\ra\geq -\ve,\,
 \forall (y,y^*)\in \partial f.
\end{equation}
Now, take
\begin{equation} \label{eq:ac4}
  \begin{array}{l}
 x_1 ^*\in \partial_{\ve_1} f(x_1),   x_2 ^* \in \partial_{\ve_1}
 f(x_2),\\
 p,q\geq 0, p+q=1,
  \end{array}
\end{equation}
and define
\begin{equation}
  \label{eq:ac5}
  \begin{array}{l}
(\bar{x},\bar{x}^*):=p ( x_1,x_1 ^*) +q (x_2, x_2 ^*),\\[.3em]
\bar{\ve}:=p \ve_1+q\ve_2+pq\la x_1 - x_2, x_1 ^* -x_2 ^*\ra .
  \end{array}
\end{equation}
Using the previous definitions, and the convexity of $f$, is trivial
to check that
\begin{equation}
  \label{eq:ac6}
   \bar{\ve}\geq 0,\; \bar{x}^*\in\partial_{\bar{\ve}} f(\bar{x}).
\end{equation}
Properties (\ref{eq:ac4},\ref{eq:ac5},\ref{eq:ac6}) will be called a
\emph{transportation formula}. If $\ve_1=\ve_2=0$, then we are using
elements in the graph of $\partial f$ to construct elements in the
graph of $\partial_{\ve} f$.
In (\ref{eq:ac5}), the product of elements in $\partial_\ve f$
appears.  This product admits the following estimation,
\begin{equation}
  \label{eq:ac7}
  x_1 ^* \in \partial_{\ve_1} f(x_1),   x_2 ^* \in \partial_{\ve_1} f(x_2)
 \Rightarrow \la x_1 - x_2, x_1 ^* - x_2 ^*\ra \geq  -(\ve_1+\ve_2).
\end{equation}
Moreover, $\partial_\ve f$ is maximal with respect to property
(\ref{eq:ac7}).  We will call property (\ref{eq:ac7}) \emph{additivity}.
The enlargement $\partial_\ve f$ can be characterized by the function
$ h_{\mathrm{FY}}$, defined in (\ref{eq:new1})
\[ x^*\in \partial_\ve f(x)\iff  h_{\mathrm{FY}}(x,x^*)\leq \la x,
x^*\ra+\ve.\]
The transportation formula (\ref{eq:ac4},\ref{eq:ac5},\ref{eq:ac6}) now
follows directly of the convexity of $ h_{\mathrm{FY}}$.  Additivity
follows from the fact that  $h_{\mathrm{FY}}\geq \J h_{\mathrm{FY}}$,
and maximality of the additivity follows from the fact that
\[
  h_{\mathrm{FY}}=\J h_{\mathrm{FY}}.
\]
Define the \emph{graph} of $\partial_\ve f$, as
\[ G(\partial_{(\cdot)} f(\cdot)):=\{(x,x^*,\ve)\,|\,
x^*\in\partial_{\ve}f(x)
\}.\]
Note that $G(\partial_{(\cdot)} f(\cdot))$ is closed.  So we say that
$\partial_\ve f$ is \emph{closed}.

Given $T:X\tos X^*$, maximal monotone, it would be desirable to
have an enlargement of $T$, say $T^\ve$, with similar properties
to the $\partial _\ve f$ enlargement of $\partial f$. With this
aim, such an object was defined in \cite{bis:97,bs:99}(in finite
dimensional spaces and in Banach spaces, respectively), for
$\ve\geq 0$,
\begin{equation}
  \label{eq:ac8}
  T^\ve(x):=\{ x^*\in X^*\,|\, \la x-y, x^*-y^*\ra\geq -\ve,\,\forall
  (y,y^*)\in T\}.
\end{equation}
The $T^\ve$ enlargement of $T$ shares with the $\partial _\ve f$
enlargement of $\partial f$ many properties: the transportation
formula, Lipschitz continuity (in the interior of its domain), and
even Br{\o}ndsted-Rockafellar property (in Reflexive Banach
spaces).  Since its introduction, it has had both theoretical and
algorithmic applications
 \cite{bis:97,bss:99a,ss99:hybrideps,ss:erbounds,rt1,rt2}.
  Even though, $T^\ve$ is \emph{not} the extension of the
construct
$\partial_{\ve }f$ to a generic maximal monotone operator. Indeed,
taking
$T=\partial f$, we obtain
\[ \partial _\ve f(x)\subseteq (\partial f)^\ve(x),\]
with examples of strict inclusion even in finite dimensional
cases~\cite{bis:97}.  Therefore, in general, $T^\ve$ lacks the
``additive'' property (\ref{eq:ac7}). The $T^\ve$ enlargement satisfy a
weaker property~\cite{bs:99}
\[
 x_1 ^* \in T^{\ve_1}(x_1),   x_2 ^* \in T^{\ve_2}(x_2)
 \Rightarrow \la x_1 - x_2, x_1 ^* - x_2 ^* \ra \geq
 -(\sqrt{\ve_1}+\sqrt{\ve_2})^2.
\]
The enlargement $T^\ve$ is also connected with a convex
function. Indeed,
\begin{eqnarray}
\nonumber
 x^*\in  T^\ve(x) & \iff & \la x-y, x^*-y^*\ra\geq -\ve,\forall
           (y,y^*)\in T\\
\nonumber
  &\iff& \sup_{(y,y^*)\in T} \la x-y, y^*-x\ra  \leq \ve.
\end{eqnarray}
Fitzpatrick function, $\varphi_T$ is the smallest element of
$\ch(T)$~\cite{fitzpatrick:1988},
and is defined as
\begin{equation}
  \label{eq:deffitz}
  \varphi_T(x,x^*):=\sup_{(y,y^*)\in T} \la x-y, y^*-x\ra+\la x,x^*\ra.
\end{equation}
Therefore,
\[
  x^*\in T^\ve(x)\iff \varphi_T(x,x^*)\leq \la x,x^*\ra+\ve.
\]
Now, the transportation formula for $T^\ve$ follows from convexity
of $\varphi_T$. In \cite{burachik-svaiter:2002} it is proven that
each enlargement $\hat{T}^\ve$ of $T$, which has a closed graph,
is nondecreasing and satisfy the transportation formula, is
characterized by a function $\hat h\in \ch(T)$, by the formula
\[ x^*\in {\hat{T}} ^\ve (x)\iff \hat h(x,x^*) \leq \la x,
  x^*\ra+\ve.\]
So, if we want to retain ``additivity'':
\[ x_1 ^* \in \hat{T}^{\ve_1}(x_1),   x_2 ^* \in \hat{T}^{\ve_2}(x_2)
 \Rightarrow \la x_1 - x_2, x_1 ^* - x_2 ^* \ra \geq
 -({\ve_1}+{\ve_2}).\]
we shall require $\hat h\geq \J \hat h$.  The enlargements in this
family, which are also maximal with respect to the additivity, are
structurally closer to the $\partial_\ve f$ enlargement, and are
characterized by $\hat h\in \ch(T)$,
\[ \hat h=\J \hat h. \]
If there were only one element in $\ch(T)$  fixed point of $\J$,
then this element would be the ``canonical'' representation of $T$
by a convex function, and the associated enlargement would be the
extension of the $\ve$-subdifferential enlargement to $T$.
Unfortunately, it is not clear whether we have uniqueness of such
fixed points.

Existence of an additive enlargement of $T$, maximal with respect
with ``additivity'' was proved in \cite{svaiter}.  The convex
representation of this enlargement turned out to be minimal in the
family $\ch_a(T)$, but the characterization of these minimal
elements of $\ch_a(T)$ as fixed point of $\J$ was lacking.

Since the function $\sigma_T$ has played a fundamental role in our
proof, we redescribe it here.  Let $\delta_T$ be the indicator
function of $T$, i.e., in $T$ its value is $0$ and elsewhere
($X\times X^*\setminus T$) its value is $+\infty$. Denote the
duality product by $\pi:X\times X^*\to \mathbb{R}$, $\pi
(x,x^*)=\la x,x^*\ra$. Then
\[ \sigma_T(x,x^*)=\mathrm{cl-conv} (\pi+\delta_T),\]
were $\mathrm{cl-conv} f$ stands for the biggest lower semicontinuous
convex function majorized by $f$.  We refer the reader to
\cite{burachik-svaiter:2002}, for a detailed analysis of this
function.

\section{Acknowledgements}

We thanks the anonymous referee for the suggestions which helped to
improve this paper.

\end{document}